\documentclass[10pt,british]{amsart}

\newif\ifpdf
\ifx\pdfoutput\undefined\pdffalse\else\pdfoutput=1\pdftrue\fi\newif\pdf
\ifpdf\relax\else
\usepackage[T1]{fontenc}
\fi

\usepackage{babel,eucal,url,amssymb,stmaryrd,booktabs,%
enumerate,amscd,paralist}

\usepackage[pagebackref]{hyperref}

\theoremstyle{plain}
\newtheorem{thm}{Theorem}[section]

\newtheorem{co}[thm]{Corollary}

\theoremstyle{definition}

\theoremstyle{remark}
\newtheorem{rem}[thm]{Remark}

\newtheorem*{ack}{Acknowledgements}

\newcommand{\Lie}[1]{\operatorname{\textsl{#1}}}

\newcommand{\GL}{\Lie{GL}}

\newcommand{\SO}{\Lie{SO}}

\newcommand{\Gtwo}{\ifmmode{{\rm G}_2}\else{${\rm G}_2$}\fi}

\newcommand{\LC}{{\nabla^g}}

\newcommand{\Nt}{\tilde\nabla}

\newcommand{\Hodge}{\mathord{\mkern1mu *}}

\date{\today}
\begin{document}

\title[$SU(3)$-instantons and $G_2, Spin(7)$-heterotic string solitons]%
{$SU(3)$-instantons and $G_2, Spin(7)$-heterotic string solitons}

\author{Petar Ivanov}
\address[P.Ivanov]{University of Sofia "St. Kl. Ohridski"\\
Faculty of Mathematics and Informatics,\\
Blvd. James Bourchier 5,\\
1164 Sofia, Bulgaria}
\email{fn10854@fmi.uni-sofia.bg}

\author{Stefan Ivanov}
\address[S.Ivanov]{University of Sofia "St. Kl. Ohridski"\\
Faculty of Mathematics and Informatics,\\
Blvd. James Bourchier 5,\\
1164 Sofia, Bulgaria}
\email{ivanovsp@fmi.uni-sofia.bg}

\begin{abstract}
Necessary and sufficient conditions to the existence of a
hermitian connection with totally skew-symmetric torsion and
holonomy contained in $SU(3)$ are given. A formula for the
Riemannian scalar curvature is obtained.  Non-compact solution to
the supergravity-type I equations of motion with
non-zero flux and non-constant dilaton is found in dimension 6.
Non-conformally flat
non-compact solutions to the supergravity-type I equations of motion with
non-zero flux and non-constant dilaton are found in dimensions
7 and 8. A Riemannian metric with holonomy contained in $G_2$
arises from our considerations and Hitchin's flow equations, which seems
to be new. Compact examples of $SU(3), G_2$ and $Spin(7)$
instanton satisfying the anomaly cancellation conditions are
presented.


\end{abstract}

\maketitle
\setcounter{tocdepth}{2}
\tableofcontents

\section{Introduction}
Supersymmetric backgrounds of string/M theory with non-vanishing fluxes are
currently an active area of study for at least two reasons. Firstly they
provide a frame work of searching for new models with realistic phenomenology
and secondly, they appear in generalizations of the AdS/CFT correspondence.

The supersymmetric geometries of the common NS-NS sector of type
IIA, IIB and heterotic/type I supergravity are analyzed in
\cite{GMW}. The bosonic geometry is of the form $\mathbb
R^{1,9-p}\times M_p$ where  the Riemannian metric $g$, the dilaton
function $\phi$ and the three form $H$ are non-trivial only on
$M_p$ but all R-R fields and fermions are set to zero in type II
theories. The type I/heterotic geometries, which is the main
object of interest in the present note, allow in addition
non-trivial gauge field $A$ with field strength $F^A$.

We recall the basic notations \cite{Str,HS,GNic,GMPW,GMW}.

We search for a solutions to lowest nontrivial order in $\alpha'$ of the equations
of motion that follow from the bosonic action
\begin{gather*}
S=\frac{1}{2k^2}\int
d^{10}x\sqrt{-g}e^{-2\phi}(Scal^g+4(\nabla^g\phi)^2-\frac{1}{12}H^2
-\alpha'Tr (F^A)^2)
\end{gather*}
which also preserves at least one supersymmetry.

The three from $H$ satisfies a modified Bianchi identity
\begin{gather}\label{modb}
dH=2\alpha'Tr(F^A\wedge F^A).
\end{gather}

The so called Bianchi identity reads
\begin{equation}\label{modb1}
dH=2\alpha'(Tr(F^A\wedge F^A)-Tr(R^{\Nt}\wedge R^{\Nt})),
\end{equation}
where $R^{\Nt}$ is the curvature of the metric connection $\Nt$
with torsion $T^{\Nt}=-H$ related to the Levi-Civita connection
$\LC$ by $\Nt=\LC-\frac12H$. The second term on the right hand
side of \eqref{modb1} is the leading string correction to the
supergravity expression arising from the anomaly cancellation but
for the consistency of the theory, a modification to the action
should be included \cite{Berg} (see also \cite{Car1},\cite{GPap}).

In terms of characteristic classes \eqref{modb1} means that $dH$
is proportional to the difference of the  first Pontrjagin 4-forms
$(p_1(A)-p_1(\Nt))$ of the connections $A,\Nt$, respectively.

A heterotic/type I geometry will preserve supersymmetry if and only if, in 10 dimensions,
there exists at least one Majorana-Weyl spinor $\epsilon$ such that the
supersymmetry variations of the fermionic fields vanish, i.e.
\begin{gather} \nonumber
\delta_{\lambda}=\nabla_m\epsilon = \left(\nabla_m^g
+\frac{1}{8}H_{mnp}\Gamma^{np} \right)\epsilon=0 \\\label{sup1}
\delta_{\Psi}=\left(\Gamma^m\partial_m\phi
+\frac{1}{12}H_{mnp}\Gamma^{mnp} \right)\epsilon=0 \\ \nonumber
\delta_{\xi}=F^A_{mn}\Gamma^{mn}\epsilon=0,
\end{gather}
where  $\lambda, \Psi, \xi$ are the gravitino, the dilatino and
the gaugino, fields, respectively.

The equations of motion corresponding to the action $S$ are
presented explicitly in \cite{GMW}. It is known \cite{Bwit,GMPW}
that the equation of motions of type I supergravity are
automatically satisfied if one imposes, in addition to the
preserving supersymmetry equations \eqref{sup1}, the modified
Bianchi identity \eqref{modb}.

According to no-go (vanishing) theorems  (a consequence of the
equations of motion \cite{FGW,Bwit}; a consequence of the
supersymmetry \cite{IP2} for SU(n)-case and \cite{GMW} for the
general case) there are no compact solutions with non-zero flux
and non-constant dilaton satisfying the supersymmetry equations
\eqref{sup1} and the modified Bianchi identity \eqref{modb}
simultaneously.

In dimensions 7 and 8 the only known heterotic/type I  solutions to
the equations of motion preserving at least one supersymmetry,
i.e. satisfying \eqref{sup1} and \eqref{modb}, are those
constructed in \cite{FNu,FN,HS} in dimension 8 and those presented
in \cite{GNic} in dimension 7. All these solutions are conformal
to a  flat space. In dimension 6, the possibility of the existence
of a non-conformally flat solution on the complex Iwasawa manifold
was discussed in \cite{Str,Car,GMW,KST}.

In the present note we concentrate our attention to find
non-compact solutions to the supergravity equations \eqref{sup1}
including the modified Bianchi identity \eqref{modb} as well as
the anomaly cancellation condition \eqref{modb1}. In dimensions 7
and 8 we find non-locally-conformally flat non-compact solutions
to the gravitino, gaugino and dilatino equations with non-zero
flux and non-constant dilaton which obey the Bianchi identities
\eqref{modb} and \eqref{modb1} and therefore satisfy the equations
of motion, due to the result in \cite{GMPW}. In dimension 6, we
present a (non-conformally-flat) non-compact solution to the
equations of motion showing that it obeys
\eqref{sup1},\eqref{modb},\eqref{modb1}. All these non-compact
solutions seem to be new.

We present a non conformally flat (resp. conformally flat) $SU(3),
G_2$ and $Spin(7)$-instantons which satisfy anomaly cancellation
condition \eqref{modb1} as well as the modified Bianchi identity
\eqref{modb}.

We obtain compact 6,7 and 8-manifolds which solve the gravitino
and gaugino equations and satisfy the compatibility conditions
\eqref{modb1} and \eqref{modb} but do not solve the dilatino
equation which is consistent with the no-go theorems.

Geometrically, the vanishing of the gravitino variation is
equivalent to the existence of a non-trivial spinor parallel with
respect to a metric connection $\nabla$ with totally skew
symmetric torsion $T=H$ which is related to the Levi-Civita
connection $\LC$ by $$ \nabla = \LC + \frac{1}{2}H. $$ The
presence of $\nabla$-parallel spinor leads to restriction of the
holonomy group $Hol(\nabla)$ of the torsion connection $\nabla$.
Namely, $Hol(\nabla)$ has to be contained in $SU(3), d=6$
\cite{Str,IP1,IP2,GIP,Car,BB,BBE}, the exceptional group $G_2,
d=7$ \cite{FI,GKMW,FI1}, the Lie group $Spin(7), d=8$
\cite{GKMW,I1}. A detailed analysis of the possible geometries is
carried out in \cite{GMW}. Complex Non-K\"ahler geometries appear
in string compactifications and are studied intensively
\cite{Str,GP,GMW,GKMW,GMPW,GP,GPap,BB,BBE}. Some types of
non-complex 6-manifold have been also invented recently in the
string theory due to the mirror symmetry and T-duality
\cite{KST,GLMW,GM,BDS,Car,Car1}.

Another special dimension is turn out to be dimension 5. The
existence of $\nabla$-parallel spinor in dimension 5 determines an
almost contact metric structure whose properties as well as
solutions to gravitino and dilatino equations are investigated in
\cite{FI,FI2}. We use these consideration in our construction in
Section~\ref{sa} of new $SU(3)$-instanton and non-compact solution to the
equations of motion in dimension 6.

Almost Hermitian manifolds with totally skew-symmetric Nijenhuis
tensor arise as target spaces of a class of (2,0)-supersymmetric
two-dimensional sigma models \cite{Pap}. For the consistency of
the theory, the Nijenhuis tensor has to be parallel with respect
to the torsion connection with holonomy contained in $SU(n)$. The
known models are those on group manifolds. We present a
6-dimensional nil-manifold as an example which is not a group
manifold.

Starting with a $SU(3)$-structure in dimension 6 we analyze the
five classes discovered recently by Chiossi and Salamon \cite{CS}
from the point of the existence of a $SU(3)$-connection having
totally skew-symmetric torsion. We obtain necessary and sufficient
conditions for the existence of a connection solving the gravitino
equation in dimension 6, i.e. the existence of a linear connection
preserving the almost hermitian structure with torsion 3-form and
holonomy contained in $SU(3)$, in terms of the given
$SU(3)$-structure and present  a formula for the Riemannian scalar
curvature (Theorem~\ref{cythm1}). It turns out that the
corresponding almost complex structure may not be integrable. In
case that the almost complex structure is integrable, we derive
that the $SU(3)$-structure is holomorphic if and only if the
corresponding hermitian structure is balanced, i.e. has co-closed
fundamental form (Corollary~\ref{colcy}). On the other hand, any
$SU(3)$-weak holonomy manifold (Nearly K\"ahler manifold)
automatically solves both the gravitino and gaugino equations. It
turns out that the Nearly K\"ahler 6-sphere $S^6$ satisfies in
addition the compatibility conditions \eqref{modb1}, \eqref{modb}.
We present a six dimensional non-conformally flat nil-manifold
$Nil^6=G/\Gamma$, i.e. a compact quotient of a nilpotent Lie group
$G$ with a discrete subgroup $\Gamma$, which  solves both the
gravitino and gaugino equations and satisfies the compatibility
conditions \eqref{modb1}, \eqref{modb} but it is neither complex
nor Nearly K\"ahler. Consequently, it does not solve the dilatino
equation \cite{Str}. However, the $SU(3)$-structure on $G$ is
half-flat and therefore it determines a Riemannian metric with
holonomy contained in $G_2$ on $G\times \mathbb R\cong \mathbb
R^7$ according to the procedure discovered by Hitchin \cite{Hit},
which seems to be a new one.

We propose a simple way to lift a 6-dimensional solution to the
gravitino and gaugino equations, i.e. a $SU(3)$-instanton solving
the gravitino equation and satisfying the conditions \eqref{modb}
(resp. \eqref{modb1}), to a $G_2$-instanton on the product with
the real line which solves the three supersymmetry equations
\eqref{sup1} as well as the compatibility condition \eqref{modb}
(resp. \eqref{modb1}). We show that $Nil^6\times\mathbb R$, (resp.
$S^6\times\mathbb R$) is a non-compact solution to the equations of motion in
dimension 7 with non-zero flux and non-constant dilaton
 which preserves at least one
supersymmetry and is not locally conformally flat (resp. locally
conformally flat). Consequently, the compact spaces $Nil^6\times
S^1, \quad S^6\times S^1$ admit a $G_2$-instanton structure
satisfying all the equations \eqref{sup1}, \eqref{modb1},
\eqref{modb} except the dilatino equation.

It turns out that any $G_2$-weak holonomy manifold
(Nearly-parallel manifold) automatically solves both the gravitino
and gaugino equations. We show that the Nearly parallel 7-sphere
satisfies in addition the compatibility conditions \eqref{modb1},
\eqref{modb}.

The same lifting procedure is applicable to the $Spin(7)$ case.
Namely, any $G_2$-instanton solving the gravitino equation and
satisfying the conditions \eqref{modb} (resp. \eqref{modb1}) can
be lifted to a $Spin(7)$-instanton on the product with the real
line which solves the three supersymmetry equations \eqref{sup1}
as well as the compatibility condition \eqref{modb} (resp.
\eqref{modb1}). We show  that $Nil^6\times\mathbb R\times\mathbb
R$, (resp. $S^6\times\mathbb R\times\mathbb R, \quad
S^7\times\mathbb R$) is a non-compact solution to the equations of
motion in dimension 8 with non-zero flux and non-constant dilaton
 which preserves at least one
supersymmetry and is not locally conformally flat (resp. locally
conformally flat). Consequently, the compact spaces $Nil^6\times
S^1\times S^1, \quad S^6\times S^1\times S^1, \quad S^7\times S^1$
admit a $Spin(7)$-instanton structure satisfying all the equations
\eqref{sup1}, \eqref{modb1}, \eqref{modb} except the dilatino
equation.

Starting with a Tanno deformed Einstein Sasaki structure in dimension 5,
we lift it to a $SU(3)$-instanton on the  product with the real line which
satisfies all the supersymmetry equations \eqref{sup1} in
dimension 6. In this way we obtain local solutions to the
equations of motion in dimension 6. Consider $S^5$ as a Sasakian
space form, i.e a Tanno deformation of the
standard Einstein-Sasaki structure, we show that $S^5\times\mathbb
R$ is a non-compact solution to the equations of motion in
dimension 6 with non-zero flux and non-constant dilaton which
preserves at least one supersymmetry. Consequently, the compact
space $S^5\times S^1$ admits a $SU(3)$-instanton structure
satisfying all the equations \eqref{sup1}, \eqref{modb1},
\eqref{modb}.

\begin{ack}
The research was done during the visit of S.I. at the Abdus Salam
International Centre for Theoretical Physics, Trieste, Italy. S.I.
thanks the Abdus Salam ICTP for providing support and an excellent
research environment. S.I. is a  member of the EDGE, Research
Training Network HPRN-CT-2000-00101, supported by the European
Human Potential Programme. The research is partially supported by
Contract MM 809/1998 with the Ministry of Science and Education of
Bulgaria, Contract 586/2002 with the University of Sofia "St. Kl.
Ohridski". We thank to Tony Pantev for his interest to this work,
for the useful  suggestions  and stimulating discussions. We are
also grateful to Jerome Gauntlett for his helpful comments and
remarks. We would like to thank the referee for his valuable
comments and suggestions on clarifying  the phenomenological
background especially  the form of the Bianchi identity including
string corrections.
\end{ack}

\section{General properties of $SU(3),G_2$ and $Spin(7)$-structures}
In this section we recall necessary properties of $SU(3),G_2$ and
$Spin(7)$ structures.
\subsection{SU(3)-structures in $d=6$}
Let $(M^6,g,J)$ be an almost Hermitian 6-manifold with Riemannian
metric $g$ and an almost complex structure $J$, i.e. $(g,J)$
define an $U(3)$-structure. The Nijenhuis tensor $N$, the K\"ahler
form $F$ and the Lee form $\theta^6$ are defined by
\begin{equation}\label{cy1}
N=[J.,J.]-[.,.]-J[J.,.]-[.,J.], \quad F=g(.,J.), \quad \theta^6(.)=\delta F(J.),
\end{equation}
respectively.

A $SU(3)$-structure is determined by an additional  non-degenerate
(3,0)-form $\Psi=\Psi^++\sqrt{-1}\Psi^-$, or equivalently by a
non-trivial spinor. To be more explicit, we may choose a local
orthonormal frame  $e_1,\dots,e_6$,  identifying it  with the dual
basis via the metric. Write $e_{i_1 i_2\dots i_p}$ for the
monomial $e_{i_1} \wedge e_{i_2} \wedge \dots \wedge e_{i_p}$. A
$SU(3)$-structure is described locally by
\begin{gather}\label{AA}
\Psi=-(e_1+\sqrt{-1}e_2)\wedge (e_3+\sqrt{-1}e_4)\wedge  (e_5+\sqrt{-1}e_6),\quad
F =-e_{12} - e_{34}- e_{56}, \\ \nonumber
\Psi^+ =-e_{135} + e_{236} + e_{146} +e_{245}, \quad
  \Psi^- =-e_{136} - e_{145} - e_{235} + e_{246}.
\end{gather}
The subgroup of $\SO(6)$ fixing the forms $F$ and  $\Psi$ simultaneously is $SU(3)$.
The two  forms $F$ and $\Psi$ determine the metric completely. The Lie algebra of $SU(3)$
is denoted $su(3)$.

The failure of the holonomy group of the Levi-Civita connection to
reduce to $SU(3)$ can be measured by the intrinsic torsion $\tau$,
which is identified with $\LC F$ or $\LC J$ and can be decomposed
into five classes \cite{CS}, $\tau \in W_1\oplus \dots \oplus
W_5$. The intrinsic torsion of an $U(n)$- structure belongs to the
first four components described by Gray-Hervella \cite{GrH}. The
five components of a $SU(3)$-structure are first described by
Chiossi-Salamon \cite{CS} (for interpretation in physics see
\cite{Car,GLMW,GM}) and are determined by $dF,d\Psi^+,d\Psi^-$ as
well as by $dF$ and $N$. We describe those of them which we will
use later.
\begin{description}
\item[$\tau \in W_1$] The class of Nearly K\"ahler (weak holonomy) manifold defined by
$dF$ to be (3,0)+(0,3)-form.
\item[$\tau \in W_2$] The class of almost K\"ahler manifolds, $dF=0$.
\item[$\tau \in W_3$] The class of balanced hermitian manifold determined by the
conditions $N=\theta^6=0$, i.e these are complex manifolds with vanishing Lee form.
These spaces are investigated in \cite{Gau1,Mi,AB}.
\item[$\tau \in W_4$] The class of locally conformally K\"ahler spaces characterized by
$dF=\theta^6\wedge F$.
\item[$\tau \in W_1\oplus W_3\oplus W_4$] The class called by Gray-Hervella $G_1$-manifolds
determined by the condition that the Nijenhuis tensor is totally skew-symmetric.
This is the precise class which we are interested in.
\end{description}
The class of a half-flat SU(3)-manifold \cite{CS} may be
characterized by the conditions $d\Psi^+=0, \quad \theta^6=0$. The
half-flat structures can be lifted to a $G_2$-holonomy metric on
the product with the real line and vice versa due to the Hitchin
theorem \cite{Hit}. In fact, many new $G_2$-holonomy metrics are
obtained in this way \cite{GLPS,BGom}.

We recall \cite{CS} that the fifth component $W_5$ and the
two scalar components of $W_1$ are determined by the expressions $\Hodge d\Psi^+\wedge\Psi^+,
\quad d\Psi^+\wedge F=W_1^+vol.,\quad
d\Psi^-\wedge F=W_1^-vol.$, respectively.

If all five components are zero then we have a Ricci-flat K\"ahler (Calabi-Yau) 3-fold.
\subsection{$G_2$-structures in $d=7$.}

Endow ${\mathbb R}^7$ with its
standard orientation and inner product.  Let $e_1,\dots,e_7$ be an oriented
orthonormal basis.
Consider the three-form $\omega$ on ${\mathbb R}^7$ given by
\begin{equation}
  \omega =e_{127} - e_{236} + e_{347}+e_{567} - e_{146} - e_{245} +
  e_{135}.\label{11}
\end{equation}
The subgroup of $\GL(7)$ fixing $\omega$ is the exceptional Lie group
$G_2$.  It is a compact, connected, simply-connected, simple Lie subgroup
of $\SO(7)$ of dimension 14~\cite{Br}. The Lie algebra is denoted by
 $g_2$ and it is isomorphic to the two forms satisfying 7 linear equations, namely
 $g_2\cong  \{\alpha\in \Lambda^2(M) \vert
  \Hodge(\alpha\wedge\omega) =- \alpha\}$
The 3-form
$\omega$ corresponds to a real spinor $\epsilon$ and therefore,
$G_2$ can be identified as the isotropy group of a non-trivial
real spinor.

The Hodge star operator supplies the 4-form $\Hodge\omega$ given by
\begin{equation}
  \Hodge\omega =  e_{3456} + e_{1457} + e_{1256}+e_{1234} + e_{2357} +
  e_{1367} - e_{2467}.\label{12}
\end{equation}
A $7$-dimensional Riemannian manifold is called a $G_2$-manifold
if its structure group reduces to the exceptional Lie group $G_2$.
The existence of a $G_2$-structure is equivalent to the existence
of a global non-degenerate three-form which can be locally written
as \eqref{11}. The 3-from $\omega$ is called  \emph{the
fundamental form} of the $G_2$-manifold \cite{Bo}. From the purely
topological point of view, a $7$-dimensional paracompact manifold
is a $G_2$-manifold if and only if it is an oriented spin
manifold~\cite{LM}.
  We will say that the pair $(M,\omega)$ is a \emph{$G_2$-manifold} with
  \emph{$G_2$-structure} (determined by) \emph{$\omega$}.

The fundamental form of a $G_2$-manifold determines a Riemannian
metric \emph{implicitly}  through
  $g_{ij}=\frac16\sum_{kl}\omega_{ikl}\omega_{jkl}$ \cite{Gr}.  This is referred to as the
  metric induced by $\omega$. We write $\LC$ for the associated Levi-Civita
  connection.

In~\cite{FG}, Fernandez and Gray divide $G_2$-manifolds into 16 classes
according to how the covariant derivative of the fundamental three-form
behaves with respect to its decomposition into $G_2$ irreducible components
(see also~\cite{CS,GKMW}).  If the fundamental form is parallel with respect to
the Levi-Civita connection, $\LC\omega=0$,
 then the Riemannian holonomy group is contained
in $G_2$. In this case the induced metric on the
$G_2$-manifold is Ricci-flat, a fact first observed by Bonan~\cite{Bo}.  It
was shown by Gray~\cite{Gr} (see also~\cite{Br,Sal}) that a $G_2$-manifold
is parallel precisely when the fundamental form is harmonic, i.e. $d\omega=d*\omega=0$.
The first
examples of complete parallel $G_2$-manifolds were constructed by Bryant and
Salamon~\cite{BS,Gibb}.  Compact examples of parallel $G_2$-manifolds were
obtained first by Joyce~\cite{J1,J2,J3} and recently by Kovalev~\cite{Kov}.

The Lee form $\theta^7$ is defined by \cite{Cabr}
\begin{equation}\label{g2li}
\theta^7=-\frac{1}{3}\Hodge(\Hodge d\omega\wedge\omega),\quad
\end{equation}
If the Lee form vanishes, $\theta^7=0$ then the $G_2$-structure is
said to be \emph{balanced}. If the Lee form is closed,
$d\theta^7=0$ then the $G_2$-structure is locally conformally
equivalent to a balanced one \cite{FI1}. If the $G_2$-structure
satisfies the condition $d\Hodge\omega=\theta^7\wedge\omega$ then
it is called \emph{integrable} and an analog of the Dolbeault
cohomology is investigated in \cite{FUg}.

\subsection{$Spin(7)$-structures in $d=8$.}

Now, let us consider ${\mathbb R}^8$ endowed with an orientation and its
standard inner product. Let $\{e_0,...,e_7\}$ be an oriented
orthonormal basis.  Consider the 4-form $\Phi$ on
${\mathbb R}^8$ given by
\begin{eqnarray}\label{1}
\Phi &=&e_{0127} - e_{0236} + e_{0347}+e_{0567} - e_{0146} - e_{0245} +
  e_{0135}
 \\ \nonumber &+&
e_{3456} + e_{1457} + e_{1256}+e_{1234} + e_{2357} +
  e_{1367} - e_{2467}.\nonumber
\end{eqnarray}
The 4-form  $\Phi$ is self-dual $*\Phi=\Phi$ and the 8-form $\Phi\wedge\Phi$ coincides with
the volume form of ${\mathbb R}^8$. The subgroup of $GL(8,R)$ which
fixes $\Phi$ is isomorphic to the double covering $Spin(7)$ of
$SO(7)$ \cite{HL}. Moreover, $Spin(7)$ is a compact
simply-connected Lie group of dimension 21 \cite{Br}. The Lie algebra of $Spin(7)$ is
denoted by $spin(7)$ and it is isomorphic to the two forms satisfying 7 linear equations, namely
 $spin(7)\cong  \{\alpha \in
\Lambda^2(M)|*(\alpha\wedge\Phi)=-\alpha\}$.

The 4-form
$\Phi$ corresponds to a real spinor $\phi$ and therefore,
$Spin(7)$ can be identified as the isotropy group of a non-trivial
real spinor.

A \emph{$Spin(7)$-structure} on an 8-manifold $M$ is by definition
a reduction of the structure group of the tangent bundle to
$Spin(7)$; we shall also say that $M$ is a \emph{$Spin(7)$
manifold}. This can be described geometrically by saying that
there exists a nowhere vanishing global differential 4-form $\Phi$
on $M$ which can be locally written as (\ref{1}). The 4-form
$\Phi$ is called the \emph{fundamental form} of the $Spin(7)$
manifold $M$ \cite{Bo}.

The fundamental form of a $Spin(7)$-manifold determines a
Riemannian metric \emph{implicitly} through
  $g_{ij}=\frac{1}{24}\sum_{klm}\Phi_{iklm}\Phi_{jklm}$ \cite{Gr}.
  This is referred to as the
  metric induced by $\Phi$.

In general, not every 8-dimensional Riemannian spin manifold $M^8$
admits a $Spin(7)$-structure. We explain the precise condition
\cite{LM}. Denote by $p_1(M), p_2(M), {\mathbb X}(M), {\mathbb
X}(S_{\pm})$ the first and the second Pontrjagin classes, the
Euler characteristic of $M$ and the Euler characteristic of the
positive and the negative spinor bundles, respectively. It is well
known \cite{LM} that a spin 8-manifold admits a $Spin(7)$
structure if and only if ${\mathbb X}(S_+)=0$ or ${\mathbb X}(S_-)=0$.
The latter conditions are equivalent to $
p_1^2(M)-4p_2(M)+ 8{\mathbb X}(M)=0$, for an appropriate choice
of the orientation \cite{LM}.

Let us recall that a $Spin(7)$ manifold $(M,g,\Phi)$ is said to be
parallel (torsion-free \cite{J2}) if the holonomy of the metric
$Hol(g)$ is a subgroup of $Spin(7)$. This is equivalent to saying
that the fundamental form $\Phi$ is parallel with respect to the
Levi-Civita connection $\nabla^g$ of the metric $g$. Moreover,
$Hol(g)\subset Spin(7)$ if and only if $d\Phi=0$ \cite{Br} (see also
\cite{Sal}) and any parallel $Spin(7)$ manifold is Ricci flat
\cite{Bo}. The first known explicit example of complete parallel $Spin(7)$
manifold with $Hol(g)=Spin(7)$ was constructed by Bryant and
Salamon \cite{BS,Gibb}.
The first compact examples of parallel $Spin(7)$ manifolds with
$Hol(g)=Spin(7)$ were constructed by Joyce\cite{J1,J2}.

There are 4-classes of $Spin(7)$ manifolds according to the
Fernandez classification \cite{F} obtained as irreducible
representations of $Spin(7)$ of the space $\nabla^g\Phi$.

The Lee form $\theta^8$ is defined by \cite{C1}
\begin{equation}\label{c2}
\theta^8 = -\frac{1}{7}*(*d\Phi\wedge\Phi)=*(\delta\Phi\wedge
\Phi).
\end{equation}
The 4 classes of Fernandez classification can be described in
terms of the Lee form as follows \cite{C1}: $W_0 : d\Phi=0; \quad
W_1 : \theta^8 =0; \quad W_2 : d\Phi = \theta^8\wedge\Phi; \quad W :
W=W_1\oplus W_2.$

A $Spin(7)$-structure of the class $W_1$ (ie
$Spin(7)$-structure with zero Lee form) is called
 \emph{a balanced $Spin(7)$-structure}.
If the Lee form is closed, $d\theta^8=0$ then the $Spin(7)$-structure is
locally conformally equivalent to a balanced one \cite{I1}.
It is shown in \cite{C1} that the Lee form of a $Spin(7)$
structure in the class $W_2$ is closed and therefore such a
manifold is locally conformally equivalent to a parallel $Spin(7)$
manifold. The compact spaces with closed but not exact Lee form
(i.e. the structure is not globally
conformally parallel) have very different topology than the parallel ones
\cite{I1}.

Coeffective cohomology and coeffective numbers of Riemannian
manifolds with $Spin(7)$-structure are studied in \cite{Ug}.

\section{The supersymmetry equations in dimensions 6, 7 and 8}

\subsubsection*{{\bf Dimension d=6.}} Necessary conditions to
have a solution to the system of dilatino and gravitino
equations in dimension 6 were derived by Strominger in \cite{Str} and then studied by
many authors \cite{GKMW,GMPW,GMW,Car,Car1,BB,BBE,BJ,GPap}

Necessary conditions to solve the gravitino equation are given in
\cite{FI}. The presence of a parallel spinor in dimension 6 leads
firstly to the reduction to $U(3)$, ie the existence of an almost
hermitian structure, secondly to the existence of a linear
connection preserving the almost hermitian structure with torsion
3-form and thirdly to the reduction of the holonomy group of the
torsion connection to SU(3).  It is shown in \cite{FI} that there
exists a unique linear connection preserving an almost hermitian
structure having totally skew-symmetric torsion  if and only if
the Nijenhuis tensor is a 3-form, i.e. the intrinsic torsion $\tau
\in W_1\oplus W_3\oplus W_4$. The torsion connection $\nabla$ is
determined by
\begin{equation}\label{cy2}
\nabla = \LC +\frac{1}{2}T, \qquad T=JdF+N=-dF(J.,J.,J.)+N=-dF^+(J.,J.,J.)+\frac{1}{4}N,
\end{equation}
where $dF^+$ denotes the (1,2)+(2,1)-part of $dF$. The (3,0)+(0,3)-part $dF^-$ is determined
completely by the Nijenhuis tensor \cite{Gau}.  If $N$ is a three form  then (see e.g.\cite{FI})
\begin{equation}\label{acy}
dF^-=-\frac{3}{4}JN.
\end{equation}
In addition, the dilatino equation forces the almost complex
structure to be integrable and the Lee form to be closed (for
applications in physics the Lee form has to be exact) determined
by the dilaton due to $\theta^6=2d\phi$ \cite{Str}.

When the almost complex structure is integrable, $N=0$, the torsion connection is also
known as the Bismut connection and was used by Bismut to prove local index theorem for
the Dolbeault operator on Hermitian non-K\"ahler manifold \cite{Bis}.
This formula was recently applied in string theory \cite{BBE}. Vanishing theorems for
the Dolbeault cohomology on compact Hermitian non-K\"ahler manifold were found in terms
of the Bismut connection \cite{AI,IP2,IP1}.
\subsubsection*{{\bf Dimension d=7.}}
The precise conditions to have a solution to the gravitino equation
in dimension 7 are found in \cite{FI}. Namely, there exists a
non-trivial parallel spinor with respect to a $G_2$-connection
with torsion 3-form $T$ if and only if there  exists a
$G_2$-structure $(\omega,g)$ satisfying the equations
\begin{equation}\label{sol7g}
d\Hodge\omega=\theta^7\wedge\Hodge\omega.
\end{equation}
In this case the torsion connection $\nabla$ is unique, the torsion 3-form $T$
is given by
\begin{equation}\label{tsol7g}
\nabla=\LC+\frac{1}{2}T, \qquad
H=T=\frac{1}{6}(d\omega,\Hodge\omega)\omega -\Hodge d\omega
+\Hodge(\theta^7\wedge\omega)
\end{equation}
and the Riemannian scalar curvature has the following expression \cite{FI1}
\begin{equation}\label{scal1}
s^g=\frac{1}{18}(d\omega,\Hodge\omega)+||\theta^7||^2 -
\frac{1}{12}||T||^2 +3\delta\theta^7.
\end{equation}
The necessary conditions to have a solution
to the system of dilatino and gravitino
equations were derived in \cite{GKMW,FI,FI1} and sufficiency was proved in \cite{FI,FI1}.
The general existence result \cite{FI,FI1} states that there exists a (local) non-trivial
solution to both dilatino and gravitino equations in dimension 7 if and only if there
there exists a $G_2$-structure $(\omega,g)$ satisfying the equations
\begin{equation}\label{sol7}
d\Hodge\omega=\theta^7\wedge\Hodge\omega, \quad
d\omega\wedge\omega=0, \quad \theta^7=2d\phi.
\end{equation}
The torsion 3-form (the flux $H$) is given by
\begin{equation}\label{tsol7}
\nabla=\LC+\frac{1}{2}T, \qquad H=T=-\Hodge d\omega +2\Hodge(d\phi\wedge\omega).
\end{equation}
The Riemannian scalar curvature satisfies
$
s^g=8||d\phi||^2 -\frac{1}{12}||T||^2 +6\delta d\phi.
$
\subsubsection*{{\bf Dimension d=8.}}
It is shown in \cite{I1} that the gravitino equation always have a
solution in dimension 8. Namely, any $Spin(7)$-structure admits a
unique $Spin(7)$-connection with totally skew-symmetric torsion
$$T=\Hodge d\Phi-\Hodge(\theta^8\wedge\Phi).$$

The necessary conditions to have a solution
to the system of dilatino and gravitino
equations were derived in \cite{GKMW,I1} and sufficiency was proved in \cite{I1}.
The general existence result \cite{I1} states that there exists a (local) non-trivial
solution to both dilatino and gravitino equations in dimension 8 if and only if there
there exists a $Spin(7)$-structure $(\Phi,g)$ with closed Lee form, $d\theta^8=0$,
which is equivalent
to the statement that the $Spin(7)$-structure is locally conformally balanced.
The torsion 3-form (the flux $H$)
and the Lee form are given by
\begin{equation}\label{tsol8}
\nabla=\LC+\frac{1}{2}T, \qquad H=T=\Hodge d\Phi- 2\Hodge(d\phi\wedge\Phi),
\quad \theta^8=\frac{12}{7}d\phi.
\end{equation}
The Riemannian scalar curvature satisfies
$
s^g=8||d\phi||^2 -\frac{1}{12}||T||^2 +6\delta d\phi.
$

In addition to these equations, the vanishing of the gaugino
variation requires the 2-form $F^A$ to be of instanton type:
(\cite{CDev,Str,HS,RC,DT,GMW})
\begin{enumerate}
\item[{\bf Case d=6}] A Donaldson-Uhlenbeck-Yau $SU(3)$-instanton i.e the gauge field $A$ is
a $SU(3)$-connection with curvature 2-form $F^A\in su(3)$.
The SU(3)-instanton condition can be written in local holomorphic coordinates
in the form \cite{CDev,Str}
\begin{equation}\label{6inst}
F^A_{\alpha\beta}=F^A_{\bar{\alpha}\bar{\beta}}=0,
F^A_{\alpha\bar\beta}F^{\alpha\bar\beta}=0.
\end{equation}
\item[{\bf Case d=7}] A $G_2$-instanton i.e. the gauge field $A$ is a
$G_2$-connection and its curvature 2-form $F^A\in g_2$. The latter can be expressed in any
of the following two equivalent ways
\begin{equation}\label{7inst}
F^A_{mn}\omega^{mn}\hspace{0mm}_p=0 \quad \Leftrightarrow
F^A_{mn}= \frac{1}{2}F^A_{pq}(\Hodge\omega)^{pq}\hspace{0mm}_{mn};
\end{equation}
\item[{\bf Case d=8}] A $Spin(7)$-instanton i.e. the gauge field $A$ is a
$Spin(7)$-connection and its curvature 2-form $F^A\in spin(7)$.
The latter is equivalent to
\begin{equation}\label{8inst}
F^A_{mn}=\frac{1}{2}F^A_{pq}\Phi^{pq}\hspace{0mm}_{mn}.
\end{equation}
\end{enumerate}

\section{Non-compact $G_2$-solution induced from a $SU(3)$-instanton}\label{sec}
In this section we show how to construct local solution to the
equations of motion in dimension 7 if we have a solution to the gravitino and gaugino
equations satisfying the modified Bianchi identity \eqref{modb} in dimension 6.

We first investigate necessary and sufficient condition to have a solution to the
gravitino equation in dimension 6, i.e. to have a $\nabla$-parallel spinors.
We prove the following
\begin{thm}\label{cythm1}
Let $(M^6,g,J,\Psi)$ be a 6-dimensional smooth manifold with a $SU(3)$-structure $(g,J,\Psi)$
or equivalently, the almost hermitian manifold $(M^6,g,J)$ has topologically trivial canonical
bundle trivialized by a (3,0)-form $\Psi$. The next two conditions are equivalent
\begin{enumerate}
\item[a)] There exists a unique $SU(3)$-connection with torsion 3-form, i.e. a
linear connection with torsion 3-form which preserves the almost hermitian structure
whose  holonomy is contained in SU(3).
\item[b)] The Nijenhuis tensor $N$ is totally-skew symmetric and the following conditions
hold
\begin{equation}\label{cycon}
d\Psi^+=\theta^6\wedge\Psi^+ -\frac{1}{4}(N,\Psi^+)\Hodge F, \quad
d\Psi^-=\theta^6\wedge\Psi^- -\frac{1}{4}(N,\Psi^-)\Hodge F.
\end{equation}
The torsion is given by
\begin{equation}\label{torcy}
T=-\Hodge dF +\Hodge(\theta^6\wedge F)+\frac{1}{4}(N,\Psi^+)\Psi^+
+ \frac{1}{4}(N,\Psi^-)\Psi^-.
\end{equation}
\end{enumerate}
The Riemannian scalar curvature is expressed in the following way
\begin{equation}\label{scal2}
s^g = \frac{1}{8}(N,\Psi^+)^2 + \frac{1}{8}(N,\Psi^-)^2 + 2
||\theta^6||^2 - \frac{1}{12}||T||^2 +3\delta\theta^6.
\end{equation}
In particular, if the structure is complex and balanced then the Riemannian scalar
curvature is non-positive.
\end{thm}
\begin{proof}
Suppose the condition a) holds. Then the Nijenhuis tensor $N$ is a three form due to Theorem 10.1
in \cite{FI}.  The conditions $\nabla\Psi=\nabla\Psi^+ =\nabla\Psi^-=0$ imply the constraints
\eqref{cycon} on the exterior derivative of the form $\Psi$. This can be checked directly
using \eqref{cy2} and \eqref{acy}.

To prove the converse we consider the product $M^7=M^6\times\mathbb R$ with the $G_2$-structure
$\omega$ defined by \cite{Hit,CS}
\begin{equation}\label{om}
\omega=-F\wedge e_7 -\Psi^+,
\end{equation}
where $e_7$ is the standard 1-form on $\mathbb R$.

We adopt the convention to
indicate  the  object on the product by
a superscript $7$, i.e. $\Hodge^7, T^7, \theta^7, \nabla^7, R^{\nabla^7}$ are the Hodge star
operator, the torsion 3-form, the Lee form, the torsion connection and its curvature, respectively.
The same objects on $M^6$ have superscript $6$.

Our idea is to check that the $G_2$-structure on the product $M^7=M^6\times\mathbb R$
satisfies the conditions \eqref{sol7g} and to apply the Friedrich-Ivanov
 result from \cite{FI} assuring the existence of $G_2$-connection $\nabla^7$
 with torsion 3-form $T^7$. We show that the torsion satisfies the condition
 $T(e_7,.,.)=0$ and therefore $\nabla^7e_7=0$. Hence, the connection $\nabla^7$ descends on
$M^6$ to a connection $\nabla^6$ which preserves the $SU(3)$-structure and has totally
skew-symmetric torsion.

We get from \eqref{om} applying \eqref{acy} and \eqref{cycon} the following sequence of
equalities
\begin{gather}\nonumber
\theta^7=-\frac{1}{3}\Hodge^7(\Hodge^7d\omega\wedge\omega)=
\frac{1}{3}\left(-\Hodge^6(\Hodge^6dF\wedge F) + \Hodge^6(\Hodge^6dF\wedge\Psi^+)e_7
-\Hodge^6(\Hodge^6d\Psi^+\wedge\Psi^+)\right)=\\ \label{thet}
= \theta^6+\frac{1}{4}(N,\Psi^-)e_7,
\end{gather}
where we used the identities
\begin{gather}\label{w2}
\Hodge^6(\Hodge^6d\Psi^+\wedge\Psi^+)=\Hodge^6(\Hodge^6(\theta\wedge\Psi^+)\wedge\Psi^+)=-2\theta,
\\ \nonumber
\Hodge^7d\omega=\Hodge^6dF-\Hodge^6d\Psi^+\wedge e_7, \quad
\theta=-\Hodge^6dF\wedge dF, \\ \nonumber \Hodge^6(\Hodge^6 dF\wedge\Psi^+)=-(\Psi^+,dF)=
-\frac{3}{4}(\Psi^-,N).\nonumber
\end{gather}
Applying the equalities $d\Hodge^6 F=-\Hodge^6
J\theta^6=\theta^6\wedge\Hodge^6 F$ and the conditions
\eqref{cycon} we obtain \eqref{sol7g}. Hence, there exists a
$G_2$-connection $\nabla^7$ with torsion 3-form $T^7$ given by
\eqref{tsol7g} on the product  $M^7=M^6\times\mathbb R$. To
compute $T^7$ we use \eqref{acy} and \eqref{thet}. We have
\begin{gather}\label{nav}
(d\omega,\Hodge^7\omega)=\Hodge^7(d\omega\wedge\omega)=2\Hodge^7(dF\wedge\Psi^+)=
-\frac{3}{2}(N,\Psi^+),
\\ \nonumber
\Hodge^7(\theta^7\wedge\omega)=\Hodge^6(\theta^6\wedge F)-\Hodge^6(\theta^6\wedge\Psi^+)e_7
+\frac{1}{4}(N,\Psi^-)\Psi^-,
\end{gather}
Now, the formula \eqref{tsol7g} and \eqref{cycon} give the desired expression \eqref{torcy}
which implies that the torsion $T^7$ does not depend on $e_7$
and therefore the connection $\nabla^7$ descends  to $M^6$.

Substituting \eqref{nav} and \eqref{thet} into \eqref{scal1} we get \eqref{scal2} for the
 Riemannian scalar curvature on the product which clearly coincides with the scalar curvature
 on $(M^6,g^6)$.
\end{proof}
\begin{co}
In dimension 6 the following conditions are equivalent:
\begin{enumerate}
\item[a)] There exists a non-trivial solution to the system of gravitino and dilatino equations
with non-zero flux $H$ and non-constant dilaton $\phi$;
\item[b)] There exists a $SU(3)$-structure $(F,\Psi)$ satisfying the conditions
\begin{equation*}
d\Psi^+=2d\phi\wedge\Psi^+, \quad d\Psi^-=2d\phi\wedge\Psi^-.
\end{equation*}
The flux $H$ is given by
\begin{equation}\label{gaun}
H=T=-\Hodge dF+2\Hodge(d\phi\wedge F).
\end{equation}
\end{enumerate}
The Riemannian scalar curvature of the solution has the expression
\begin{equation}\label{scal}
s^g=8||d\phi||^2 -\frac{1}{12}||T||^2 +6\delta d\phi.
\end{equation}
\end{co}
A part of the necessary conditions we presented are known \cite{Str,GMW,GKMW,GMPW,Car}. The formula
\eqref{gaun} was discovered in \cite{GKMW}, the first formula in \eqref{w2} has already
appeared in \cite{Car}.
\begin{co}\label{colcy}
A closed $SU(3)$-structure ($d\Psi^+=d\Psi^-=0$) admits a
$SU(3)$-connection with torsion 3-form if and only if the
corresponding almost Hermitian structure is a balanced Hermitian
structure, ($N=\theta^6=0$).

In particular,  a holomorphic $SU(3)$-structure on a complex manifold
supports a linear connection with torsion 3-form and holonomy contained in $SU(3)$ if and
only if it is balanced. In the latter case the Riemannian scalar curvature is non-positive.

Consequently, such a structure is half-flat and therefore it determines a Riemannian metric
with  holonomy contained in $G_2$  on the product with the real line .
\end{co}
\begin{rem}
We note the coincidence of the formulas for the Riemannian scalar curvature of solutions
to the gravitino and dilatino equations in dimensions 6,7, and 8.
Actually, it is proved in \cite{IP1} (see also \cite{GPap}) that the SU(n)-geometry
arising from any solution to the gravitino and dilatino equations satisfies the identity
\begin{equation}\label{eqms1}
Ric^g(X,Y)-\frac{1}{4}\sum_{i,j}H(X,e_i,e_j)H(Y,e_i,e_j)+2\LC_X\LC_Y\phi +\frac{1}{4}\sum_idH(X,JY,e_i,Je_i)=0
\end{equation}
which is consistent with the first equation of motion. The trace in \eqref{eqms1} gives \eqref{scal}
due to the  identity
\begin{equation*}
 \frac{1}{4}\sum_idH(X,JY,e_i,Je_i)=8||d\phi||^2+4\delta d\phi-\frac{1}{3}||H||^2
\end{equation*}
shown in \cite{AI} (see also (3.24) in \cite{IP1}).
\end{rem}
\begin{rem}
We note that the Riemannian scalar curvature for a half-flat $SU(3)$-structure is
computed in \cite{GLMW}. On the other hand, not every half-flat structure admits $SU(3)$-connection
with torsion 3-form
since it may have nonzero $W_2$ component. For example, the $SU(3)$-structure on the
nilpotent  Lie algebras described in \cite{CS}, Example 2 and 3, are half-flat but do not
admit $SU(3)$ connection with torsion 3-form since $d\Psi^-\not= fF\wedge F$.
\end{rem}
Another consequence of Theorem~\ref{cythm1} is the following
\begin{thm}\label{cyth2}
Let $(M^6,g,J,F,\Psi)$ be a smooth almost complex 6-manifold with totally skew-symmetric
Nijenhuis tensor equipped with a $SU(3)$-structure $\Psi$ solving the gravitino equation,
i.e. the conditions \eqref{cycon} hold. Assume also the conditions
\begin{equation}
d\theta^6=0, \quad (N,\Psi^+)=0, \quad (N,\Psi^-)= const.\not=0.
\end{equation}
\begin{enumerate}
\item[i)] Then the $G_2$-structure $\omega=-F\wedge e_7 -\Psi+$ defined on the product
$M^7=M^6\times\mathbb R$ solves both the gravitino and dilatino equations with non-constant
dilaton.
\item[ii)] If in addition the torsion connection $\nabla^6$ is a $SU(3)$-instanton then
the corresponding torsion connection $\nabla^7$ is a $G_2$-instanton.
\item[iii)]
Suppose moreover that the torsion connection $\nabla^6$ satisfies
the modified Bianchi identity \eqref{modb}, (resp. \eqref{modb1})
with $F^A=R^{\nabla^6}$ . Then $\nabla^7$ also obeys \eqref{modb},
(resp. \eqref{modb1}) with $F^A=R^{\nabla^7}$ and therefore solves
the  equations of motion with non zero flux and  non-constant
dilaton provided $\theta^6$ is exact.
\end{enumerate}
\end{thm}
\begin{proof}
Equations \eqref{thet} and \eqref{nav} imply \eqref{sol7} due to
the conditions of the theorem. We know from Theorem~\ref{cythm1}
that $T^7=T^6,\quad R^{\nabla^7}=R^{\nabla^6},\quad
R^{\Nt^7}=R^{\Nt^6}$ and consequently, $Tr(R^{\nabla^7}\wedge
R^{\nabla^7})=Tr(R^{\nabla^6}\wedge R^{\nabla^6}), \quad
 Tr(R^{\Nt^7}\wedge R^{\Nt^7})=Tr(R^{\Nt^6}\wedge
R^{\Nt^6})$. A glance at the structure of Lie algebras $su(3)$ and
$g_2$ implies that $R^{\nabla^7}$ satisfies the $G_2$- instanton
equations \eqref{7inst} provided $R^{\nabla^6}$ obeys the
$SU(3)$-instanton equations \eqref{6inst}
\end{proof}

\section{Non-compact $Spin(7)$-solution induced from a $G_2$-instanton}\label{sec1}
In this section we shall show that a $G_2$-instanton on $(N^7,\omega)$ induces a
$Spin(7)$-instanton on the product $N^8=N^7\times \mathbb R$.

We denote the Hodge star operator on $N^7$ by $\Hodge^7$.
On the product $N^8=N^7\times \mathbb R$ there exists a $Spin(7)$-structure defined by
\begin{equation}\label{sg1}
\Phi=e_0\wedge\omega + *^7\omega,
\end{equation}
where $e_0=dt$ is the standard 1-form on $\mathbb R$. We indicate  the  object on the product by
a superscript $8$, i.e. $\Hodge^8, T^8, \theta^8, \nabla^8, R^{\nabla^8}$ are the Hodge star
operator, the torsion 3-form, the Lee form, the torsion connection and its curvature, respectively.
The same objects on $N^7$ have superscript $7$.
\begin{thm}\label{th1}
Suppose $(N^7,\omega^7, g^7, \nabla^7, T^7)$ is a smooth $G_2$-manifold which solves the
gravitino equation, i.e. \eqref{sol7g} holds.
 Then the $Spin(7)$-structure $\Phi$ on the product
 $N^7\times \mathbb R$ determined with \eqref{sg1} has the properties
\begin{equation*}
\theta^8=\frac{6}{7}\theta^7+\frac{1}{7}(d\omega,\Hodge^7\omega)e_0,
\quad T^8=T^7, \quad R^{\nabla^8}=R^{\nabla^7},\quad
Tr(R^{\nabla^8}\wedge R^{\nabla^8})=Tr(R^{\nabla^7}\wedge
R^{\nabla^7}).
\end{equation*}
Assume in addition the conditions
\begin{equation*}
d\theta^7=0, \qquad (d\omega,\Hodge^7\omega)=const.\not=0.
\end{equation*}
\begin{enumerate}
\item[i)] Then the $Spin(7)$-structure $\Phi$ defined on the product
$N^8=N^7\times\mathbb R$ solves both the gravitino and dilatino equations with non-zero
flux and  non-constant dilaton.
\item[ii)] If in addition the torsion connection $\nabla^7$ is a $G_2$-instanton then
the corresponding torsion connection $\nabla^8$ is a $Spin(7)$-instanton.
\item[iii)]
Suppose moreover that the torsion connection $\nabla^7$ satisfies
the modified Bianchi identity \eqref{modb}, (resp. \eqref{modb1})
with $F^A=R^{\nabla^7}$. Then $\nabla^8$ also obeys \eqref{modb},
(resp. \eqref{modb1}) with $F^A=R^{\nabla^8}$  and therefore
solves the equations of motion with non-zero flux and non-constant
dilaton provided $\theta^7$ is exact.
\end{enumerate}
\end{thm}
\begin{proof}
Take the exterior derivative in \eqref{sg1} and use \eqref{sol7} to get
 the identity $d\Phi=-e_0\wedge d\omega +\theta^7\wedge\Hodge^7\omega$.
The latter yields
\begin{gather}\label{lee}
\Hodge^8d\Phi=-\Hodge^7d\omega - e_0\wedge\Hodge^7(\theta^7\wedge\Hodge^7\omega)\\ \label{lee1}
\Hodge^8(\theta^7\wedge\Phi)=-\Hodge^7(\theta^7\wedge\omega)-e_0\wedge\Hodge^7(\theta^7\wedge
\Hodge^7\omega)\\ \label{lee2}
\theta^8=-\frac{1}{7}\Hodge^8\left(\Hodge^8d\Phi\wedge\Phi\right)=\\ \nonumber
-\frac{1}{7}\left[\Hodge^8\left(e_0\wedge\Hodge^7 d\omega\wedge\omega\right)-
e_0\wedge\Hodge^7\left(\theta^7\wedge\Hodge^7\omega\right)\wedge\Hodge^7\omega
-\Hodge^8(\Hodge^7d\omega\wedge\Hodge^7\omega)\right]=\\ \nonumber
-\frac{1}{7}\left[\Hodge^7\left(\Hodge^7 d\omega\wedge\omega\right)-
\Hodge^7\left(\Hodge^7\left(\theta^7\wedge\Hodge^7\omega\right)\wedge\Hodge^7\omega\right)
-(d\omega,\Hodge^7\omega)e_0\right]
=\frac{6}{7}\theta^7+\frac{1}{7}(d\omega,\Hodge^7\omega)e_0,\nonumber
\end{gather}
where we used the conditions \eqref{sol7}, \eqref{g2li} and the general identity
$\Hodge^7\left(\Hodge^7\left(\gamma\wedge\Hodge^7\omega\right)\wedge\Hodge^7\omega\right)=3\gamma$
valid for any 1-form $\gamma$ on $(M^7,\omega)$.

Substitute \eqref{lee}, \eqref{lee1} and \eqref{lee2} into the formula \eqref{tsol8} and compare
the result with \eqref{tsol7} to get $T^8=T^7$.

The vector field $e_0$ is parallel with respect to the Levi-Civita
connection of $g^8$ and satisfies $T^8(e_0,.,.)=T^7(e_0,.,.)=0$.
Therefore $\nabla^8e_0=0$ yielding
$R^{\nabla^8}=R^{\nabla^7},\quad R^{\Nt^8}=R^{\Nt^7}$.
Consequently,
$$\sum_0^7R^{\nabla^8}_{mnkl}\Phi^{mn}\hspace{0.0mm}_{ij}=
\sum_1^7R^{\nabla^7}_{mnkl}(\Hodge^7\omega)^{mn}\hspace{0.0mm}_{ij}=
2R^{\nabla^7}_{ijkl}=2R^{\nabla^8}_{ijkl},$$ since $R^{\nabla^7}$
is a $G_2$-instanton. Clearly, the modified Bianchi identity
\eqref{modb}, (resp. \eqref{modb1}) is satisfied for
$F^A=R^{\nabla^8}$.
\end{proof}
The reverse procedure to find types of $SU(3)$-structures on 6-manifold induced by
different types of $G_2$-structures on 7-manifold is discussed recently in \cite{BJ,BDS,GLMW,GM}.

\section{Examples}\label{examp}
Theorem~\ref{cyth2} and Theorem~\ref{th1} allow us to produce a
number of examples of $G_2$ and $Spin(7)$-instantons and solutions
to the equations of motion with gauge connection $A=\nabla$
 starting from certain types of almost complex 6-manifold
or certain types of $G_2$ manifolds.

We recall the well known curvature identity
\begin{equation}\label{bas1}
R^{\nabla}(X,Y,Z,V)=R^{\Nt}(Z,V,X,Y) + \frac12dT(X,Y,Z,V).
\end{equation}
It  helps us to handle the Bianchi identity  \eqref{modb1} with
$F^A=R^{\nabla}$ provided the next equality  holds
\begin{equation}\label{bas2}
R^{\nabla}(X,Y,Z,V)=R^{\nabla}(Z,V,X,Y).
\end{equation}
Combine \eqref{bas1} and \eqref{bas2} to get
\begin{equation}\label{bas3}R^{\nabla}=R^{\Nt}+\frac12dT.
\end{equation}
In view of \eqref{bas3}, the Bianchi identity \eqref{modb1} with
$A=\nabla$ takes the form
\begin{equation*} dT=\alpha'\left(Tr(R^{\nabla}\wedge dT) - \frac12Tr(dT\wedge dT)\right).
\end{equation*}

Clearly, the  condition \eqref{bas2} is sufficient a $SU(3)$
(resp. $G_2, Spin(7)$) -connection to satisfy the $SU(3)$ (resp.
$G_2,Spin(7)$) -instanton condition \eqref{6inst} (resp.
\eqref{7inst},\eqref{8inst}). The symmetry \eqref{bas2} of the
curvature of a metric connection $\nabla$ with torsion 3-form $T$
holds exactly when $\nabla T$ is a 4-form which is equivalent to
the condition $\LC T=\frac{1}{4}dT$ \cite{I2}. In particular,
if the torsion is $\nabla$-parallel, $\nabla T=0$, then we have
the additional relations
\begin{gather}\label{partor}
R^{\LC}_{ijkl}=R^{\nabla}_{ijkl}-\frac{1}{2}T_{ijm}T_{kl}\hspace{0mm}^m
-\frac{1}{4}T_{jkm}T_{il}\hspace{0mm}^m
-\frac{1}{4}T_{kim}T_{jl}\hspace{0mm}^m;\\\nonumber
dT_{ijkl}=2\left(T_{ijm}T_{kl}\hspace{0.mm}^m+
T_{jkm}T_{il}\hspace{0.mm}^m +
T_{kim}T_{jl}\hspace{0.mm}^m\right).
\end{gather}

\subsection{$(SU(3), G_2, Spin(7))$-instanton and conformally flat non-compact solution}
Any Nearly-K\"ahler 6-manifold is an $SU(3)$-instanton since the
torsion $T=\frac{1}{4}N=J\LC J$ is $\nabla$- parallel \cite{Kir},
(see also \cite{BM,FI}) and therefore the curvature $R^{\nabla}$
satisfies \eqref{bas2}.

Take $\Psi^+=dF$ we obtain a $SU(3)$-instanton solving the gravitino and gaugino equations
according to Theorem~\ref{cyth2} which, however, does not solve the dilatino equation
since the almost complex structure is not integrable \cite{Str}.

There are known only four compact Nearly K\"ahler 6-manifolds, namely $S^6$, $S^3\times S^3$,
$\mathbb CP^3$  and the flag
$F=U(3)/(U(1)\times U(1)\times U(1))$ \cite{Hit,Sal1}.

We consider the six-sphere $(S^6,J,g)$ endowed with  the standard
nearly K\"ahler structure $(g,J)$ inherited from the imaginary
octonions in $\mathbb R^7$ \cite{Gr}. We claim that
$(S^6,g,J,\nabla,A=\nabla)$ satisfies both the modified Bianchi
identity \eqref{modb} and the anomaly cancellation condition
\eqref{modb1}.

It is well known that any 6-dimensional Nearly K\"ahler manifold is Einstein and
of constant type. Consequently, the following identities hold \cite{FI}
\begin{gather*}
T_{ijm}T_{kl}\hspace{0.mm}^m=a^2\frac{1}{2}(g_{ik}g_{jl}-g_{jk}g_{il}
-F_{ik}F_{jl}+F_{jk}F_{il}), \quad dT=a^2F\wedge F=-2a^2\Hodge F,
\end{gather*}
where $a^2$ is a non-zero constant which can be identified with
the Riemannian scalar curvature $s^g$, $15a^2=s^g$. Applying the
fact that $(S^6,g)$ is a space of constant sectional curvature,
i.e. $R^{\LC}_{ijkl}=\frac{1}{2}a^2(g_{jk}g_{il}-g_{ik}g_{jl})$
and \eqref{bas3}, we calculate the Pontrjagin forms
\begin{gather*}
16\pi^2p_1(\nabla)= Tr(R^{\nabla}\wedge
R^{\nabla})=R^{\nabla}_{ijab}R^{\nabla}_{kl}\hspace{0.mm}^{ab}dx^i\wedge
dx^j\wedge dx^k\wedge dx^l = -\frac{3a^2}{4} dT;\\
16\pi^2p_1(\Nt)= Tr(R^{\Nt}\wedge R^{\Nt})=\frac{9a^2}{4} dT.
\end{gather*}
\begin{rem}\label{rem1}
Observe that if we rescale the metric homothetically by a constant $c$,
  $\bar g=e^{2c}g$ then the new torsion $\bar T=e^{2c}T$ in the case of $SU(3)$-structure and
$\bar T=e^{4c}T$ in the case of $G_2$ or $Spin(7)$-structure while the
  Pontrjagin 4-forms remain unchanged, $Tr(R^{\bar{\nabla}}\wedge R^{\bar{\nabla}})=
  Tr(R^{\nabla}\wedge R^{\nabla}),\quad Tr(R^{\bar{\Nt}}\wedge R^{\bar{\Nt}})=
  Tr(R^{\Nt}\wedge R^{\Nt})$ (see \cite{Car1} for more precise discussion of this phenomena).
  Hence, if $dT$ is proportional to the difference of
  the Pontrjagin 4-forms with a constant then we can always rescale the structure by a suitable
  constant in order to get the formulas \eqref{modb} and \eqref{modb1}.
\end{rem}
Keeping Remark~\ref{rem1} in mind we obtain
\begin{thm}
The Nearly K\"ahler 6-sphere solves the gravitino equation, the
gaugino equation with $F^A=R^{\nabla}$  and satisfies the modified
Bianchi identity \eqref{modb} and \eqref{modb1} with negative
$\alpha'$. Consequently,
\begin{enumerate}
\item[a)] the product $(S^6\times \mathbb R, \omega, A=\nabla^7)$ with the
$G_2$-structure described in section~\ref{sec} solves all the supersymmetry equations
\eqref{sup1} with non-zero flux, non-constant dilaton and satisfies the
Bianchi identity \eqref{modb}, \eqref{modb1}.
Therefore it solves the equations of motion in dimension 7.

The product $(S^6\times S^1, \omega, A=\nabla^7)$ is a compact
space solving locally the supersymmetry equations \eqref{sup1}
which satisfies the
 Bianchi identity \eqref{modb}, \eqref{modb1} in dimension 7.
\item[b)]
the product $(S^6\times \mathbb R\times \mathbb R, \Phi,
A=\nabla^8)$ with the $Spin(7)$-structure described in
section~\ref{sec1} solves the supersymmetry equations \eqref{sup1}
with non-zero flux, non-constant dilaton and satisfies the Bianchi
identity \eqref{modb}, \eqref{modb1}. Therefore it solves the
equations of motion  in dimension 8.

The product $(S^6\times S^1\times S^1, \Phi, A=\nabla^8)$ is a
compact space solving locally the supersymmetry equations
\eqref{sup1} which satisfies the
 Bianchi identity \eqref{modb}, \eqref{modb1} in dimension 8.
\end{enumerate}
\end{thm}
Similarly to the SU(3)-case, any $G_2$-weak holonomy manifold (nearly-parallel $G_2$ manifold)
is a $G_2$-instanton. Indeed, it is well known that any 7-dimensional
nearly-parallel $G_2$ manifold is Einstein and the following identities hold
\begin{equation}\label{g2par}
d\omega=-\lambda\Hodge\omega, \quad (d\omega,\Hodge\omega)=-\lambda,
\end{equation}
where $\lambda^2=\frac{8}{21}s^g$ is a non-zero constant. The
torsion $T=-\frac{1}{6}\lambda\omega$ is $\nabla$-parallel
\cite{FI}. Hence, any nearly-parallel $G_2$-manifold is a
$G_2$-instanton which solves   the gravitino equation and the
gaugino equation with $A=\nabla^7$ but does not solve the dilatino
equation
 according to the result in \cite{FI1}.

There are many known examples of compact nearly parallel $G_2$-manifolds:   $S^7$,
$SO(5)/SO(3)$ \cite{BS,Sal}, the Aloff-Wallach spaces $N(g,l)=SU(3)/U(1)_{gl}$ \cite{CMS}, any
Einstein-Sasakian and any 3-Sasakian 7-manifold \cite{FK,FKMS}.

We consider the seven sphere $(S^7,\omega,g)$ endowed with the
standard nearly-parallel $G_2$-structure induced by the octonions
in $\mathbb R^8$, namely, consider the seven sphere as a totally
umbilical hypersurface in $\mathbb R^8$ \cite{FG}. Clearly
$(S^7,\omega,g,\nabla, A=\nabla)$ is a $G_2$-instanton. We claim
that it satisfies the modified Bianchi identity \eqref{modb} and
the anomaly cancellation condition \eqref{modb1}. Indeed, we
easily calculate from \eqref{partor} applying \eqref{g2par},
\eqref{bas3}, the fact that $(S^7,g)$ is a space of  constant
sectional curvature and some $G_2$-algebra, that
\begin{gather*}
16\pi^2p_1(\nabla)= Tr(R^{\nabla}\wedge
R^{\nabla})=R^{\nabla}_{ijab}R^{\nabla}_{kl}\hspace{0.mm}^{ab}dx^i\wedge
dx^j\wedge dx^k\wedge dx^l =
-\frac{\lambda^4}{8.27}\Hodge\omega=-\frac{\lambda^2}{36} dT;\\
16\pi^2p_1(\Nt)= Tr(R^{\Nt}\wedge R^{\Nt}) =
\frac{1}{54}\lambda^4\Hodge\omega= \frac19\lambda^2dT.
\end{gather*}
We obtain using  Remark~\ref{rem1}
the following
\begin{thm}
The nearly-parallel 7-sphere solves the gravitino equation, the
gaugino equation with $F^A=R^{\nabla}$ and satisfies the modified
Bianchi identity \eqref{modb} and \eqref{modb1} with negative
$\alpha'$.

Consequently, the product $(S^7\times \mathbb R, \Phi,
A=\nabla^8)$ with the $Spin(7)$-structure described in
section~\ref{sec1} solves the supersymmetry equations \eqref{sup1}
with non-zero flux, non-constant dilaton and satisfies the Bianchi
identity \eqref{modb}, \eqref{modb1}. Therefore it solves the
equations of motion  in dimension 8.

The product $(S^7\times S^1, \Phi, A=\nabla^8)$ is a compact space
solving locally the supersymmetry equations \eqref{sup1} which
satisfies the
 Bianchi identity \eqref{modb}, \eqref{modb1} in dimension 8.
\end{thm}
We note that these sphere-solutions are (locally) conformally flat.

\subsection{$(SU(3), G_2, Spin(7))$-instanton and non-conformally flat non-compact solution}
In this section we present a non-locally-conformally flat solution starting with
a nilpotent 6-dimensional Lie group.

Let $G$ be the six-dimensional connected simply connected and nilpotent
 Lie group, determined by the left invariant 1-forms $\{e_1,\dots,e_6\}$
 such that
\begin{gather}\label{in1}
de_2=de_3=de_6=0,\\\nonumber
de_1=e_3\wedge e_6,\quad
de_4= e_2\wedge e_6, \quad
de_5= e_2\wedge e_3.
\end{gather}
In terms of the standard coordinates $x_1,\dots,x_6$ on $\mathbb R^6$ the left invariant forms
$\{e_1,\dots,e_6\}$ are described by the expressions
\begin{gather}\label{his1}
e_2=dx_2, e_3=dx_3, e_6=dx_6 \\\nonumber
e_1=dx_1 - x_6dx_3,\quad
e_4=dx_4 -x_6dx_2, \quad
e_5=dx_5+x_2dx_3.
\end{gather}
Consider the metric on $G\cong \mathbb R^6$ defined by
$g=\sum_{i=1}^6e_i^2$, or equivalently
\begin{gather}\label{hmet}
ds^2=dx_1^2+(1+x_6^2)dx_2^2+(1+x_2^2+x_6^2)dx_3^2+dx_4^2+dx_5^2+dx_6^2 \\
-x_6(dx_1dx_3+dx_2dx_4)+x_2dx_3dx_5. \nonumber
\end{gather}
Let $(F,\Psi)$ be the $SU(3)$-structure on $G$ given by \eqref{AA}. Then$(G,F,\Psi)$ is
an almost complex manifold with a $SU(3)$-structure.

We show below that
this space is a new non-conformally flat $SU(3)$-instanton
solving both the gravitino and gaugino equations satisfying the modified Bianchi identity
\eqref{modb} as well as the anomaly cancellation condition \eqref{modb1}
 but not solving the dilatino equation.

We compute the Riemannian curvature $R^g$. The general Koszul formula
\begin{gather}\label{kz}
2g(\LC_XY,Z)=Xg(Y,Z)+Yg(X,Z)-Zg(X,Y)\\ \nonumber +g([XY],Z)-g([Y,Z],X)-g([X,Z],Y)
\end{gather}
gives the following essential non-zero terms
\begin{gather} \nonumber
2\LC_{e_6}e_3=e_1, \quad  2\LC_{e_2}e_3=-e_5, \quad 2\LC_{e_6}e_2=e_4, \\ \label{levc}
2\LC_{e_3}e_6=-e_1, \quad 2\LC_{e_3}e_2=e_5, \quad 2\LC_{e_2}e_6=-e_4, \\ \nonumber
2\LC_{e_1}e_6=-e_3, \quad  2\LC_{e_5}e_2=e_3, \quad 2\LC_{e_4}e_6=-e_2, \\ \nonumber
\end{gather}
Then we obtain $R^g(e_5,e_6,e_2,e_1)=-\frac{1}{4}\not=0$. Hence,
the metric is not locally conformally flat since the Weyl tensor
$W^g(e_5,e_6,e_2,e_1)=R^g(e_5,e_6,e_2,e_1)=-\frac{1}{4}\not=0$.

It is easy to verify using  \eqref{cy1} and \eqref{in1} that
\begin{gather}\label{in2}
dF=-3e_{236}, \quad  N=-\Psi^-, \quad d\Psi^-=\Hodge F, \quad
(N,\Psi^-)=-4, \\ \nonumber \theta^6=d\Psi^+=(N,\Psi^+)=0.
\end{gather}

Hence, $(G,\Psi,g,J)$ is neither complex nor Nearly K\"ahler manifold but it fulfills
the conditions \eqref{cycon} of Theorem~\ref{cythm1} and
therefore there exists a $SU(3)$-holonomy connection with torsion 3-form on $(G,\Psi,g,J)$.
The expression
\eqref{torcy} and \eqref{in2} give
\begin{equation}\label{tor}
T=-2e_{145}+e_{136}+e_{235}-e_{246}, \quad dT=-2(e_{1256}+e_{3456}+e_{1234})=2\Hodge F.
\end{equation}
Plug  \eqref{levc} and \eqref{tor} into \eqref{cy2} to get that the nonzero essential terms
of the torsion connection are
\begin{gather}\label{tor1}
\nabla_{e_1}e_6=-e_3, \qquad  \nabla_{e_5}e_2= e_3, \qquad \nabla_{e_4}e_6=-e_2, \\\nonumber
\nabla_{e_4}e_5=-e_1, \qquad \nabla_{e_5}e_1=-e_4, \qquad \nabla_{e_1}e_4=-e_5.
\end{gather}
It follows from \eqref{tor1} and \eqref{tor} that the torsion tensor $T$ as well as the
Nijenhuis tensor $N$ are parallel with
respect to the connection $\nabla$.  Hence, $\nabla$ defines an $SU(3)$-instanton.

To verify the Bianchi identities  for $H$ we calculate the curvature $R^{\nabla}$ by
means of \eqref{tor1}. We obtain the following non-zero terms
\begin{gather}\nonumber
R^{\nabla}(e_6,e_2,e_6,e_2)=R^{\nabla}(e_6,e_3,e_6,e_3)=
R^{\nabla}(e_2,e_3,e_2,e_3)=1\\\label{7curv}
R^{\nabla}(e_4,e_5,e_4,e_5)=R^{\nabla}(e_4,e_1,e_4,e_1)=R^{\nabla}(e_5,e_1,e_5,e_1)=1\\
\nonumber
 R^{\nabla}(e_2,e_6,e_5,e_1)=R^{\nabla}(e_3,e_6,e_4,e_5)=R^{\nabla}(e_2,e_3,e_1,e_4)=-1.
\end{gather}

Applying \eqref{7curv}, \eqref{tor} and \eqref{bas3} it is
straightforward to compute the first Pontrjagin 4-forms
$p_1(\nabla)$ and $p_1(\Nt)$. Compare the result with the second
equality in \eqref{tor} to get
\begin{equation}\label{pont}
dT=\frac12Tr(R^{\nabla}\wedge R^{\nabla})=- Tr(R^{\Nt}\wedge
R^{\Nt}).
\end{equation}

The coefficient of the structure equations of the Lie algebra given by \eqref{in1}
are integers. Therefore, the well-known theorem of Malcev \cite{Mal} states that
the group $G$ has a uniform discrete subgroup $\Gamma$ such that $Nil^6=G/\Gamma$ is a
compact 6-dimensional nil-manifold. The $SU(3)$-structure, described above, descends to $Nil^6$
and therefore we obtain a compact $SU(3)$-instanton.
With the help of Remark~\ref{rem1} and \eqref{pont} we derive
from Theorem~\ref{cyth2} and Theorem~\ref{th1}
the following
\begin{thm}
The non conformally flat almost hermitian 6-manifold
$(G,g,F,\Psi,A=\nabla)$ solves the gravitino and gaugino equations
and satisfies the Bianchi identity \eqref{modb}, \eqref{modb1}
with positive $\alpha'$. Consequently,
\begin{enumerate}
\item[a)] the product $(G\times \mathbb R, \omega, A=\nabla^7)$ with the
$G_2$-structure described in section~\ref{sec} solves all the supersymmetry equations
\eqref{sup1} with non-zero flux, non-constant dilaton and satisfies the
Bianchi identity \eqref{modb}, \eqref{modb1}.
Therefore it solves the equations of motion in dimension 7.

The product $(Nil^6\times S^1, \omega, A=\nabla^7)$ is a compact
space solving locally the supersymmetry equations \eqref{sup1}
which satisfies the Bianchi identity \eqref{modb} and
\eqref{modb1} in dimension 7;
\item[b)]
the product $(Nil^6\times \mathbb R\times \mathbb R, \Phi,
A=\nabla^8)$ with the $Spin(7)$-structure described in
section~\ref{sec1} solves the supersymmetry equations \eqref{sup1}
with non-zero flux, non-constant dilaton and satisfies the Bianchi
identity \eqref{modb}, \eqref{modb1}. Therefore it solves the
equations of motion  in dimension 8.

The product $(Nil^6\times S^1\times S^1, \Phi, A=\nabla^8)$ is a
compact space solving locally the supersymmetry equations
\eqref{sup1} which satisfies the
 Bianchi identity \eqref{modb}, \eqref{modb1} in dimension 8.
\end{enumerate}
\end{thm}
The $G_2$-analog of the Dolbeault cohomology on $G_2$-manifold was
studied on $Nil^6\times S^1$ in \cite{FUg}.
\begin{rem}
The space $(G,g,J)$ is an example of an almost complex 6-manifold
with totally skew-symmetric Nijenhuis tensor $N$ and zero Lee form
$\theta^6$ which is neither complex nor Nearly K\"ahler but it is
half-flat and therefore it determines a Riemannian metric with
holonomy contained in $G_2$ on $G\times \mathbb R\cong\mathbb R^7$
\cite{Hit} which seems to be new. The explicit expression of this
metric can be found solving the Hitchin flow equations
$$dF=\frac{\partial(\Psi^+)}{\partial t}, \quad d\Psi^-=-F\wedge
\frac{\partial(F)}{\partial t},$$ where the $SU(3)$-structure
depends on a real parameter $t \in \mathbb R$ \cite{Hit}.
$G_2$-holonomy metrics arising from types of Hermitian 6-manifolds
are studied recently in \cite{AS}.
\end{rem}
\begin{rem}
The torsion tensor $T$ as well as the Nijenhuis tensor $N$ of $(Nil^6,g,J)$ are parallel
with respect to the torsion connection $\nabla$ but $\nabla R^{\nabla}\not=0$ since
the space is not naturally reductive due to the inequality
$-1=g([e_3,e_6],e_1)\not=-g([e_3,e_1],e_6)=0$. Hence, $(Nil^6,g,J)$ is an example
of a compact non-naturally reductive almost Hermitian 6-manifold with totally
skew-symmetric Nijenhuis tensor which is neither complex nor Nearly K\"ahler.
The torsion as well as the Nijenhuis tensor are parallel
with respect to the torsion connection.
\end{rem}
\begin{rem}
Spaces for which the covariant derivative of the torsion is a four
form (in particular zero) become automatically of 'instanton
type'. Spaces with parallel torsion are studied in \cite{Ag1} in
connection with string model; almost Hermitian 6-manifolds with
parallel torsion are investigated very recently (after the first
version of the present article was posted to the arXiv) in
\cite{AFr}.
\end{rem}
\begin{rem}
In general, on any almost hermitian manifold the Nijenhuis tensor
$N$ is the (3,0)+(0,3)-part of the torsion of any linear
connection $\nabla$ compatible with the almost hermitian structure
(see e.g. \cite{Gau}). Therefore, the condition $\nabla T=0$
always implies $\nabla N=0$ because $\nabla$ preserves the type
decomposition induced from the almost complex structure.
\end{rem}
\section{Almost contact metric structures and non-compact $SU(3)$-solutions in dimension 6}\label{sa}
We construct in this section a new non-compact solution to the type I-supergravity equations of motion
in dimension 6. We derive our solution from  Sasakian structures in dimension 5.

Solutions to the gravitino and dilatino equations in dimension 5 are investigated
in \cite{FI,FI2}. In dimension five any solution to the gravitino
equation, i.e. any parallel spinor with respect to a metric connection with
torsion 3-form defines an almost contact metric structure $(g,\xi,\eta,\psi)$ which
is preserved by the torsion connection \cite{FI,FI2}.
It is shown in \cite{FI2} that solutions to the both gravitino and
dilatino equations are connected with a special type 'conformal' transformations
of the almost contact structure introduced in \cite{FI2}.

We recall that an almost contact metric structure consists of
an odd dimensional manifold $M^{2k+1}$ equipped with a Riemannian metric $g$, vector field $\xi$
of length one, its dual 1-form $\eta$ as well as an endomorphism $\psi$ of the tangent bundle
such that
\begin{equation*}
\psi(\xi)=0, \quad \psi^2=-id +\eta\otimes\xi, \quad g(\psi.,\psi.)=g(.,.)-\eta\otimes\eta.
\end{equation*}
The Nijenhuis tensor $N$ and the fundamental form $F$ of an almost contact metric structure
are defined by
\begin{equation*}
F(.,.)=g(.,\psi.), \quad N=[\psi,\psi]+d\eta\otimes\xi.
\end{equation*}
There are many special types of almost contact metric structures.
We introduce those which are relevant to our considerations:
\begin{enumerate}
\item[-] normal almost contact structures determined by the condition $N=0$;
\item[-] contact metric structures characterized by $d\eta=2F$;
\item[-] quasi-Sasaki structures, $N=0, dF=0$.
Consequently, $\xi$ is a Killing vector \cite{Bl};
\item[-] Sasaki structures, $N=0, d\eta=2F$.
Consequently, $\xi$ is a Killing vector \cite{Bl}.
\end{enumerate}
An almost contact metric structure admits a linear connection
$\nabla$ with torsion 3-form preserving the structure, i.e.
$\nabla g=\nabla\xi=\nabla\psi=0$, if and only if the Nijenhuis
tensor is totally skew-symmetric and the vector field $\xi$ is
Killing vector field \cite{FI}. In this case the torsion
connection is unique. The torsion $T$ of $\nabla$ on a Sasakian
manifold is expressed by $T=\eta\wedge d\eta=2\eta\wedge F$ and
the torsion $T$ is $\nabla$-parallel, $\nabla T=0$ \cite{FI}.

We restrict our attention to the Sasakian manifold in dimension
five.

The spinor bundle $\Sigma$ of a 5-dimensional contact metric spin
manifold decomposes under the action of the fundamental 2-form
$F^5$ into the sum $\Sigma=\Sigma^0\oplus\Sigma^1\oplus\Sigma^2,
dim \Sigma^0=dim\Sigma^2=1, dim\Sigma^1=2$. Spinors of type
$\Sigma^1$ parallel with respect to the torsion connection on
quasi-Sasakian 5-manifold are studied in \cite{FI2}.

We are interested in $\nabla^5$-parallel spinors of type $\Sigma^0
{\rm or} \Sigma^2$ on Sasakian 5-manifold. We recall (\cite{FI1},
Theorem 9.2) that a 5-dimensional simply connected Sasakian
manifold admits $\nabla^5$-parallel spinor of type $\Sigma^0 {\rm
or} \Sigma^2$ if and only if the Riemannian Ricci tensor $Ric^g$
has the form
\begin{equation}\label{ric1}
Ric^g=6g-2\eta\otimes\eta.
\end{equation}
The Tanno deformation of a Sasakian structure satisfying
\eqref{ric1}, defined by the formulas
\begin{equation*}
\psi'=\psi, \quad \xi'=\frac34\xi, \quad \eta'=\frac43\eta,\quad
g'=\frac43g+\frac49\eta\otimes\eta,
\end{equation*}
yields an Einstein-Sasakian structure with Ricci tensor
$Ric^{g'}=4g'$ and vice versa.

We may choose locally an orthonormal basis
$e_1,e_2,e_3,e_4,e_5=\xi$ such that
\begin{equation}\label{sas1}
F^5=e_1\wedge e_2 + e_3\wedge e_4, \quad d\eta=de_5=2F^5, \quad
T^5=2\eta\wedge F^5= 2e_5\wedge(e_1\wedge e_2 + e_3\wedge e_4).
\end{equation}
Let $M^6=M^5\times\mathbb R$. We indicate the objects on $M^6$
 with a superscript 6. Consider the
product Riemannian manifold $M^6=M^5\times\mathbb R$
 with the product metric and the
compatible almost complex structure $J$ determined by the
fundamental form
\begin{equation}\label{sas2}
F^6=F^5 + e_5\wedge e_6,
\end{equation}
where $e_6=dt$ is the standard 1-form on $\mathbb R$.

The identity \eqref{sas2} yields $dF^6=2F^5\wedge e_6=2e_6\wedge
F^6$. This equality tell us that the almost hermitian manifold
$(M^6,g^6,F^6)$ is locally conformally K\"ahler. In particular,
the almost complex structure is integrable and the Lee form
$\theta^6=2e_6=2dt$ is a closed 1-form on $\mathbb R$. It is easy
to check that the torsion $T^6$ of the corresponding Bismut
connection $\nabla^6$ is determined by the equality $T^6=T^5$.
Consequently, $\nabla^6 T^6=0$.  The Bismut connection defines an
$SU(3)$-instanton on $M^5\times\mathbb R$ which solves the three
supersymmetry equations \eqref{sup1} provided $M^5$ is Sasakian
5-manifold whose Riemannian Ricci tensor satisfies  \eqref{ric1}.

\subsection{Non-conformally flat local $SU(3)$-solutions on $S^5\times S^1$.}
Let $S^5$ be  the five-sphere  with the standard Einstein Sasakian
structure induced on $S^5$ by the usual complex structure on
$\mathbb C^3$ considering $S^5$ as a totally geodesic hypersurface
in the complex space $\mathbb C^3$. We consider
$(S^5,g,\psi,\eta,\xi)$ as a Sasakian space form, i.e. a Tanno
deformation of the standard Einstein-Sasakian structure on $S^5$.
We may assume that the Riemannian curvature is given by (see e.g.
\cite{Bl}
\begin{gather}\label{sas}
R^g_{ijkl}=\frac43\left(g_{jk}g_{il}-g_{ik}g_{jl}\right)+
\frac13\left(F_{kj}F_{li}-F_{ki}F_{lj}+2F_{ij}F_{lk}\right)+\\\nonumber
\frac13\left(\eta_i\eta_kg_{jl} - \eta_j\eta_kg_{il} +
\eta_j\eta_lg_{ik} - \eta_i\eta_lg_{jk}\right).
\end{gather}
Consider the locally conformally K\"ahler structure on
$S^5\times\mathbb R$ determined by \eqref{sas2}. We claim that the
corresponding Bismut connection satisfies the modified Bianchi
identity \eqref{modb} as well as the anomaly cancellation
condition \eqref{modb1}. Indeed, the equation \eqref{sas1} yields
$dT^6=dT^5=d\eta\wedge d\eta = 4F^5\wedge F^5$.

Use \eqref{partor},  apply \eqref{sas1} and
\eqref{sas}, to get
\begin{gather*}
R^{\nabla^6}_{ijkl}=\frac43\left(g_{jk}g_{il}-g_{ik}g_{jl} +
\eta_i\eta_kg_{jl} - \eta_j\eta_kg_{il} + \eta_j\eta_lg_{ik} -
\eta_i\eta_lg_{jk}\right) +\frac16dT^6_{ijkl}.
\end{gather*}
The latter equality as well as \eqref{bas3} help  to calculate the
Pontrjagin forms
\begin{equation*}
16\pi^2p_1(\nabla^6)= Tr(R^{\nabla^6}\wedge R^{\nabla^6})
=-\frac{8}{3}dT^6, \quad 16\pi^2p_1(\Nt^6)=\frac{16}{3}dT^6.
\end{equation*}

Keeping Remark~\ref{rem1} in mind we obtain
\begin{thm}
The Sasakian space form $(S^5,g^5,\psi,\eta,\xi)$  solves the
gravitino  equation and satisfies the modified Bianchi identity
\eqref{modb} and \eqref{modb1} for $F^A=R^{\nabla^5}$ with
negative $\alpha'$. Consequently,
\begin{enumerate}
\item[a)] the product $(S^5\times \mathbb R, g^6, F^6,\nabla^6,A=\nabla^6)$ solves all the supersymmetry equations
\eqref{sup1} with non-zero flux, non-constant dilaton and satisfies the
Bianchi identity \eqref{modb}, \eqref{modb1}.
Therefore it solves the equations of motion in dimension 6.

The product $(S^5\times S^1, g^6,F^6,A=\nabla^6)$ is a compact
space solving locally the supersymmetry equations \eqref{sup1}
which satisfies the
 Bianchi identity \eqref{modb}, \eqref{modb1} in dimension 6.

\item[b)] the product $(S^5\times \mathbb R\times\mathbb R,\omega,A=\nabla^7)$ with the
$G_2$-structure described in section~\ref{sec} solves all the supersymmetry equations
\eqref{sup1} with non-zero flux, non-constant dilaton and satisfies the
Bianchi identity \eqref{modb}, \eqref{modb1}.
Therefore it solves the equations of motion in dimension 7.

The product $(S^5\times S^1\times S^1,\omega,A=\nabla^7)$ is a
compact space solving locally the supersymmetry equations
\eqref{sup1} which satisfies the
 Bianchi identity \eqref{modb}, \eqref{modb1} in dimension 7.
\item[c)]
the product $(S^5\times \mathbb R\times \mathbb R\times\mathbb
R,\Phi,A=\nabla^8)$ with the $Spin(7)$-structure described in
section~\ref{sec1} solves the supersymmetry equations \eqref{sup1}
with non-zero flux, non-constant dilaton and satisfies the Bianchi
identity \eqref{modb}, \eqref{modb1}. Therefore it solves the
equations of motion  in dimension 8.

The product $(S^5\times S^1\times S^1\times S^1,\Phi,A=\nabla^8)$
is a compact space solving locally the supersymmetry equations
\eqref{sup1} which satisfies the
 Bianchi identity \eqref{modb}, \eqref{modb1} in dimension 8.
\end{enumerate}
\end{thm}
\begin{rem}
Note that all compact examples we have presented in Section~\ref{examp} and
Section~\ref{sa} solve the
supersymmetry equations only locally since the closed Lee form $\theta$ is actually
a closed 1-form on a circle and therefore it can not be exact. This is consistent with
the  vanishing results  claiming that there are no compact solutions
with globally defined non-constant dilaton and non-zero flux in type II and type I
supergravities.
\end{rem}

\bibliographystyle{hamsplain}

\providecommand{\bysame}{\leavevmode\hbox to3em{\hrulefill}\thinspace}






\end{document}